\documentstyle[12pt,righttag]{amsart}
\makeatletter

\newtheorem{prop}{Proposition}[section]
\newtheorem{thm}[prop]{Theorem}

\newtheorem{lem}[prop]{Lemma}

\theoremstyle{definition}
\newtheorem{defn}{Definition}[section]

\theoremstyle{remark}
  
\newcommand{\bbD}{{\Bbb{D}}}
\newcommand{\bbN}{{\Bbb{N}}}
\newcommand{\bbR}{{\Bbb{R}}}
\newcommand{\bbT}{{\Bbb{T}}}
\newcommand{\bbZ}{{\Bbb{Z}}}

\newcommand{\calB}{{\cal{B}}}
\newcommand{\calG}{{\cal{G}}}

\newcommand{\calD}{{\cal{D}}}

\newcommand{\calP}{{\cal{P}}}

\newcommand{\calQ}{{\cal{Q}}}

\newcommand{\frakB}{{\frak{B}}}
\newcommand{\frakA}{{\frak{A}}}

\newcommand{\lam}{\lambda}
\newcommand{\sig}{\sigma}

\newcommand{\eps}{\epsilon}

\newcommand{\Ome}{\Omega}

\newcommand{\ome}{\omega}
\newcommand{\del}{\delta}

\newcommand{\vphi}{\varphi}

\newcommand{\Del}{\Delta}

\renewcommand{\mod}{\operatorname{mod}}
\renewcommand{\Im}{\operatorname{Im}}

\newcommand{\mes}{\operatorname{mes}}

\newcommand{\Ker}{\operatorname{Ker}}
\newcommand{\diag}{\operatorname{diag}}

\renewcommand{\max}{\operatornamewithlimits{max}}

\renewcommand{\Re}{\operatorname{Re}}

\newcommand{\vl}[1]{{\bf #1}}

\newcommand{\tp}[1]{``{\it #1\/}''}

\newcommand{\lb}{\label}

\newcommand{\ti}{\tilde}

\newcommand{\bs}{\backslash}

\newcommand{\emp}{\emptyset}

\newcommand{\bi}[1]{\bibitem[#1]{#1}}

\makeatletter
\def\@currentlabel{2.1}\label{e:dispaa}
\def\@currentlabel{2.21}\label{e:dispau}
\def\@currentlabel{2.22}\label{e:dispav}
\def\@currentlabel{2.23}\label{e:dispaw}
\def\@currentlabel{2.24}\label{e:dispax}

\makeatother

\makeatletter
\def\alphenumi{%
  \def\theenumi{\alph{enumi}}%
  \def\p@enumi{\theenumi}%
  \def\labelenumi{(\@alph\c@enumi)}}
\makeatother

\begin{document}
\begin{center}

{\huge Evolutionary Semigroups and Lyapunov Theorems in Banach Spaces}

\bigskip\bigskip
{\large Yuri Latushkin}

\bigskip
{\large  University of Missouri, Columbia, Missouri 65211}

\bigskip
{\large and}

\bigskip
{\large Stephen Montgomery-Smith}
\footnote{This author was supported by
  the National Science
  Foundation under the grant  DMS-9201357.}

\bigskip
{\large  University of Missouri, Columbia, Missouri 65211}

\end{center}
\baselineskip=16pt

\begin{abstract}
We present a spectral mapping theorem for continuous semigroups of
operators on any Banach space $E$.  The condition for the
hyperbolicity of a semigroup on $E$ is given in terms of
the generator of an evolutionary semigroup acting in the space of
$E$-valued functions.  The evolutionary semigroup generated by
the propagator of a nonautonomous differential
equation in $E$ is also studied.  A ``discrete'' technique for the
investigating of the evolutionary semigroup is developed and applied to
describe the hyperbolicity (exponential dichotomy) of the nonautonomous
equation.
\end{abstract}

\section{Introduction}
Let us consider an autonomous differential equation $v' =Av$ in a Banach
space $E$, where $A$ is a generator of $C_0$-semigroup $\{e^{tA}\}_{t\geq
0}$.  Denote, as usual, $s(A) =\sup \{\Re \lam: \lam \in \sig (A)\}$
and $\ome (A) =\inf \{ \ome : \|e^{tA}\| \leq Me^{\ome t}\}$.

A classical result of A.~M.~Lyapunov (see, e.g., \cite{9}) shows
that for any bounded operator $A\in B(E)$\
the spectrum $\sig (A)$ of $A$ is
responsible for the asymptotic behavior of the solution
$y(t)=e^{At}y(0)$ of the above equation.  For example, if $\sigma (A)$
is contained in the left half-plane, that is
$s(A) <0$,
then the trivial solution is uniformly
asymptotically stable, that is
$\ome (A) <0$,
and $\|e^{tA}\| \to 0$ as $t\to \infty$.  This fact
follows from the spectral mapping theorem (see, e.g., \cite{21}):
\begin{equation}
\sig (e^{tA}) \bs \{0\} =\exp t \sig (A), \quad t>0,
\lb{smt}
\end{equation}
which always holds for bounded $A$.

For unbounded $A$, equation \eqref{smt} is not always true.  Moreover,
there are examples of generators $A$ (see \cite{21}) such that even
$s(A) <0$ does not guarantee
$\ome(A)<0$ and $\|e^{tA}\|\to 0$ as $t\to\infty$.  Since $\sig (A)$ does
not characterize the asymptotic behavior of the solution
$v(\cdot)$, we would
like to find some other characterization that still does not involve
solving the differential equation.

In this paper we solve precisely this problem in the following manner.
Consider the space $L_p (\bbR;E)$ of $E$-valued functions for $1\leq p
< \infty$, or the space $C_0 (\bbR; E)$ of continuous vanishing at $\pm
\infty$ functions on $\bbR$, and the semigroup $\{e^{tB}\}_{t\geq 0}$ of
evolutionary operators
\begin{equation}
(e^{tB}f) (x) =e^{tA} f(x-t), t \geq 0,
\lb{evop}
\end{equation}
generated by the operator $B$ that is the closure of $-\frac{d}{dx} +A$,
$x\in\bbR$.  It turns out that it is $\sig (B)$ in $L_p (\bbR;E)$ (or in
$C_0(\bbR; E)$) that is responsible for the asymptotic behavior of $v(\cdot)$
in $E$.  For example, $s(B) =\ome(A)$, and $s(B) <0$ in $L_p (\bbR;E)$ or
$C_0 (\bbR;E)$ implies $\|e^{tA}\| \to 0$ as $t\to\infty$ on $E$.

The order of the proofs is as follows. First, we consider the
evolutionary semigroup $\{e^{tB}\}$ in the space $L_p([0,2\pi];E)$
or $C([0,2\pi];E)$ of $2\pi$-periodic functions. We prove that
$1\notin \sigma (e^{2\pi A})$ in $E$ is equivalent to
$1\notin \sigma (e^{2\pi B})$
or $0\notin \sigma (B)$ in $L_p([0,2\pi];E)$ or
$C([0,2\pi];E)$. The main part of the proof
uses a modification of an idea due to C.~Chicone and R.~Swanson
\cite{6}. Next, using this result, we give a variant of the
spectral mapping theorem for a semigroup $\{e^{tA}\}$ in a Banach
space $E$. This spectral mapping theorem is a direct generalization
of L.~Gearhart's spectral mapping theorem for Hilbert spaces
(see, e.g.,
\cite[p. 95]{21}),
and is related to the spectral mapping theorem of
G.~Greiner \cite[p. 94]{21}. Finally, using a simple change
of variables arguments, we prove that
$\sigma (e^{tA})\cap \bbT =\emptyset$, $t\neq 0$
in $E$, $\bbT =\{z: |z|=1\}$
is equivalent to
$\sigma (e^{tB})\cap \bbT =\emptyset$, $t\neq 0$
which in turn is equivalent to $0\notin\sigma (B)$
in $L_p(\bbR; E)$ or $C_0(\bbR; E)$.

We will also consider the well-posed nonautonomous equation $v'=A(x) v$
in $E$, and its associated
evolutionary family $\{U(x,s)\}_{x \geq s}$,
which can be thought
as a propagator of this equation, that is
$v(x) =U(x,s)v(s)$. We assume that
$U(\cdot\, , \, \cdot)$ is strongly continuous and satisfies the usual
(\cite[p.~89]{31}) algebraic properties of the propagator.
Instead of the semigroup given by \eqref{evop} we consider on $L_p
(\bbR;E)$ or $C_0 (\bbR; E)$ the evolutionary semigroup
\begin{equation}
(e^{tG} f)(x) =U(x,x-t) f(x-t), \quad x\in\bbR, t\geq 0.
\lb{eve}
\end{equation}
We will show that $\sig (G)$ characterizes the asymptotic behavior of
$v(\cdot)$ and prove the spectral mapping theorem for the semigroup
$\{e^{tG}\}$.

We will be considering not only stability but also the exponential
dichotomy (hyperbolicity) for the solutions of the equation $v'=A(x) v$ or
evolutionary family $\{U(x,s)\}$.  We say, that
an evolutionary family $\{U(x,s)\}_{x\geq s}$ is
(spectrally) hyperbolic if there exists a continuous in the strong sense,
bounded, projection-valued function $P:\bbR \to B(E)$ such that: a) The
norm of the restrictions $U(x,s) \big| \Im P(s)$ (resp.\ $U(x,s) \big|
\Ker P(s)$) exponentially decreases (resp.\ increases) with $x-s$, and b)
$\Im U(x,s) \big| \Ker P(s)$ is dense in $\Ker P(x)$. Note, that b)
automatically follows from a) if the operators $U(x,s)$ are invertible
and defined for all $(x,s) \in\bbR^2$.
This happens, in particular, if $U(\cdot\, , \, \cdot)$ is a
norm-continuous propagator for the differential
equation $v'=A(x)v$ with continuous and bounded
$A:\bbR\to B(E)$. For this case the hyperbolicity of
$\{U(x,s)\}$ coincides with the exponential dichotomy (see, e.g.,
\cite{9}) of the equation.
However, generally a) does not imply b) (see \cite{27}).

Exponential dichotomy in the theory of differential equations with bounded
coefficients on $E$ is a major tool used for proving instability theorems for
nonlinear equations, and for showing
existence and uniqueness of bounded solutions and
Green's functions, etc.\ (see, e.g. \cite{7,9}).  The spectral
mapping theorem for the semigroup \eqref{eve} which is given here allows
one to extend these ideas to the case of unbounded coefficients.

It turns out that the spectrum $\sig (e^{tG})$ for nonperiodic $A(\cdot)$
plays the same role in the description of exponential dichotomy as the
spectrum of the monodromy operator does in the usual Floquet theory for
the periodic
case.  That is, the condition $\sig (e^{tG}) \cap \bbT =\emptyset$,
$t\neq 0$, or equivalently $0\notin \sig (G)$, is
equivalent to the (spectral) hyperbolicity of the evolutionary family
$\{U(x,s)\}_{x\geq s}$.

Showing that the hyperbolicity
$\sig (e^G) \cap\bbT =\emptyset$ of the operator $e^G$ implies
the hyperbolicity of $\{U(x,s)\}$\ is a delicate matter.
It turns out that the Riesz projection $\calP$
for $e^G$
on $L_p(\bbR,E)$ or $C_0 (\bbR;E)$,
that corresponds to $\sig (e^G) \cap
\bbD$, $\bbD=\{z: |z|< 1 \}$,
has the form $(\calP f) (x) =P(x)f(x)$.  Here $P(\cdot)$ is a continuous
in the strong sense,
projection-valued function, that defines the hyperbolicity of $\{U(x,s)\}$.
Note that $e^G =aR$, where $a$ is an operator of multiplication by
the function $a(x) =U(x,x-1)$, that is $(af)(x)=a(x)f(x)$,
and $R$
is a translation operator
$(Rf)(x) =f(x-1)$, $x\in\bbR$.
Therefore, $e^G$ falls into the class of so-called
weighted translation operators, which are well-understood in the case that
$E$\ is
Hilbert space and $p=2$ (see \cite{1,2,18,27}, and
also \cite{8,23} and references therein). If
$U(\cdot\, , \, \cdot)$ is norm-continuous,
then $\calP$ is an operator from
a $C^*$-algebra generated by $R$ and the $C^*$-algebra $\frakA_{nc}$
of operators of multiplication by the norm-continuous, bounded functions
from $\bbR$
to $B(E)$.  The techniques from the theory of weighted translation
operators (see \cite{1,2,18,27}) allows one to conclude
that $\calP  \in \frakA_{nc}$.
This technique is not applicable to the case where $\{U(\cdot,\cdot)\}$
is only strongly-continuous, nor also to the case when $E$ is not
Hilbert space.

In this paper we
present some new approaches, which allows one to derive the above
result for any Banach space $E$ and
is new even for the Hilbert space case and
when it is only known that $U(\cdot , \cdot)$ is
strongly continuous.
The main idea is to
``discretize'' the operator $aR$, that is to represent it by the family of
operators $\pi_x (a) S$, $x\in\bbR$, acting on the
``discrete" space $l_p (\bbZ;E)$.
Here $S: (v_n)_{n\in\bbZ} \mapsto (v_{n-1})_{n\in\bbZ}$ is the shift
operator and $\pi_x(a) :(v_n)_{n\in\bbZ}\mapsto (a(x+n)
v_n)_{n\in\bbZ}$ is a diagonal operator on $l_p (\bbZ;E)$.
This idea goes back to the theory of
regular representations of $C^*$-algebras
\cite{26}, and is related to works
\cite{1,2,13,16,17,18}.
As a result we prove that
$\sig (aR) \cap \bbT =\emptyset$ in $L_p(\bbR;E)$ implies $\sig (\pi_x
(a) S)\cap \bbT =\emptyset$ in $l_p (\bbZ; E)$ for each $x\in\bbR$, and
derive from this fact that
$\calP\in \frakA$, where $\frakA$ is the set of
bounded functions  $a: \bbR\to B(E)$ which
are continuous in strong operator
topology on $B(E)$.

We point out that the investigation of evolutionary operators
\eqref{evop}--\eqref{eve} has a long history, probably starting from
\cite{14} (see also \cite{10,11,19,22}).
Recently significant progress has been made in the papers
\cite{3,4,25,27}.  It is these papers
that essentially motivated and
influenced this present work.

Finally, the results of this article can be generalized to the case of the
variational equation $v'(t) =A(\vphi^t x)v(t)$ for a flow $\{\vphi^t\}$
on a compact metric space $X$, or to the linear skew-product flow $\hat
\vphi^t: X \times E \to X \times E: (x,v) \mapsto (\vphi^t x, \Phi
(x,t)v)$, $t\geq 0$ (see \cite{6,12,18,29,30} and references contained
therein).  Here $\Phi : X \times \bbR_+ \to L(E)$ is a cocycle over
$\vphi^t$, that is, $\Phi (x,t+s)=\Phi (\vphi^t x, s) \Phi (x,t)$.  Let
us recall (see \cite{29,30}) that part of the purpose
of the theory of linear skew-product flows was to
be able to handle the equation $v' =A(t) v$ in the case
when $A(\cdot)$ is almost-periodic .

To answer the question when $\hat
\vphi^t$ is hyperbolic (or Anosov), instead of \eqref{eve} one considers
the semigroup of so called weighted composition operators (see
\cite{6,15,18}) on $L_p(X;\mu;E)$:
\begin{equation}
(T^t f)(x) =\Big( \dfrac{d\mu \circ \vphi^{-t}}{d\mu}\Big)^{1/p} \Phi
(\vphi^{-t}x,t)f(\vphi^{-t} x), \quad x\in X,\; t \geq 0.
\lb{4}
\end{equation}
Here $\mu$ is a $\vphi^t$-quasi-invariant Borel measure on $X$.  As
above, the condition $\sig (T^t) \cap \bbT =\emptyset$ is equivalent to
the spectral hyperbolicity of the linear skew-product flow $\hat
\vphi^t$.  The spectral hyperbolicity coincides with the usual
hyperbolicity if $\Phi (x,t)$, $x\in X$, $t\geq 0$ are invertible or
compact operators.  A detailed investigation of weighted composition
operators and their connections with the spectral theory of linear
skew-product flows and other questions of dynamical system theory may be
found in \cite{18} (see also \cite{27}).

We will use the following notations: $\bbD=\{z:|z| < 1\}$;
$\bbT=\{z:|z|=1\}$;
``$\Big|$" denotes the restriction of an operator;
$\calD(\cdot)=\calD_F(\cdot)$ denotes the domain of
an operator in a space $F$; $\sigma(\cdot)=\sigma(\cdot\, ; F)$
denotes the spectrum;
$\sigma_{ap}(\cdot)= \sigma_{ap}(\cdot\, ; F)$
denotes the approximative point spectrum;
$\sigma_r(\cdot)=\sigma_r(\cdot\, ; F)$ denotes
the residual spectrum;
and $\rho(\cdot)=\rho(\cdot\, ; F)$ denotes the resolvent
set of an operator on $F$.
For an operator $A$ in $E$ we denote by $\cal A$
the operator of multiplication by $A$ in a space of $E$-valued
functions: $({\cal A}f)(x)=Af(x),~f: \bbR\to E$.

The authors would very much like to thank C.~Chicone for help and
suggestions, and R.~Rau for many illuminating discussions.

The authors also would like to thank the referee for
the suggestion to shorten the proof of Theorem 2.5 and Remark below.

\section{Autonomous Case}

Let $A$ be a generator of a $C_0$-semigroup $\{e^{tA}\}_{t\geq 0}$ on a
Banach space $E$.  The semigroup is called {\it hyperbolic} if $\sig
(e^{tA}) \cap \bbT =\emptyset$ for $t\neq 0$.  In this section we will
characterize the hyperbolicity of the semigroup $\{e^{tA}\}$ in terms
of evolutionary semigroup $\{e^{tB}\}_{t \geq 0}$ and its generator
$B$.  This semigroup acts by the rule $(e^{tB} f)(x) =e^{tA} f(x-t)$ on
functions $f$ with values in $E$.  In Subsection~2.1 we consider
$\{e^{tB}\}$ acting on the space $L_p ([ 0, 2\pi]; E)$ and $C([0,2\pi];E)$.
In Subsection~2.3 $\{e^{tB}\}$ acts on $L_p (\bbR;E)$ and
$C_0(\bbR;E)$.  Subsection~2.2 is devoted to a spectral mapping
theorem for $\{e^{tA}\}_{t\geq 0}$ on $E$ which generalizes the spectral
mapping theorem of L.~Gearhart for Hilbert space.

\section*{2.1. Periodic Case}\lb{sub2.1}

Let $F$ denote one of the spaces
$L_p ([0,2\pi],E)$, $1\leq p < \infty$
or  $C([0,2\pi],E)$
of $2\pi$-periodic $E$-valued functions $f$, $f(0)=f(2\pi)$.
Consider the evolutionary semigroup
$\{e^{tB}\}_{t\geq 0}$ acting on $F$, defined by the rule
$$
(e^{tB} f) (x) =e^{tA}
f([x-t] (\mod 2\pi)),~~x\in [0,2\pi].
$$
Of course,
$[0,2\pi]$ here was chosen for convenience, and for a semigroup $(e^{tB}
f)(x) =e^{tA} f$ $([x-t] (\mod t_0))$
the proofs below remain the same for any
$t_0 > 0$.

Note that $e^{tB}$ in $F$ is a product of two commuting
semigroups $(U^t f) (x) =f([x-t] (\mod 2\pi))$ and $(e^{t{\cal A}}f)(x)
=e^{tA}f(x)$.  Hence
the generator $B$ is the closure of the operator
\begin{equation}
(B_0f)(x) =-\dfrac{d}{dx} f(x) +Af(x),
\lb{genB}
\end{equation}
where $B_0$ is defined on the core $\calD (B_0)$ of $B$ (see
\cite[p.~24]{21}).  Moreover, $\calD_F(B_0) =\calD_F (-d/ dx) \cap
\calD_F ({\cal A})$, where the derivative $d/dx$ is taken in the strong sense in
$E$, and $\calD_F (B_0)=\{f: [0,2\pi] \to \calD (A) \big| f\in F \mbox{ is
absolutely continuous}, \frac{d}{dx} f \in F,
\mbox{ and } {\cal A}f\in F\}$.

Since $Be^{ik\cdot} f(\cdot) =e^{ik\cdot} (B-ik)f(\cdot)$,
$k\in\bbZ$, for the operator $(L_kf)(x)=e^{ikx}f(x)$
one has $BL_k=L_k(B-ik)$.
Therefore,
the spectrum
$\sig (B)$ in $F$
is invariant under translations by $i$.

We will need the following Lemma.

\begin{lem}\lb{l2.1}
If $1 \in \sig_{ap} (e^{2\pi A})$ in $E$, then $0 \in \sig_{ap} (B)$ in
$F$.
\end{lem}

\begin{pf}
Fix $m\in\bbN$, $m\geq 2$.  Since $1\in \sig_{ap} (e^{2\pi A})$,
we can choose $v\in E$ such that $\|v\|_E =1$ and $\|v-e^{2\pi A}v\|_E <
\frac1m$.  Note also that $\|e^{2\pi A} v\|_E\geq 1-\frac1m$.

Let $\alpha: [0,2\pi]\to [0,1]$ be any smooth function with bounded
derivative such that $\alpha (x)=0$ for $x\in [0,\frac{2\pi}{3}]$ and
$\alpha(x)=1$ for $x\in [\frac{4\pi}{3}, 2\pi]$.
Define a function $g: [0,2\pi]\to E$ by the the formula
\begin{equation}
g(x) =[1-\alpha(x)]
e^{(2\pi +x)A} v+\alpha (x) e^{xA}v,~~x\in [0,2\pi].
\lb{2.1}
\end{equation}
Note that $g(0) =g(2\pi) =e^{2 \pi A} v$.  Obviously, $g\in F$.
Also,
$$
(e^{tB} g)(x) = [1-\alpha (x-t)]e^{(2\pi +x)A} v+\alpha (x-t) e^{xA} v,$$
$g\in \calD_F(B),$
and
\begin{align}
(Bg) (x) &= \alpha' (x) e^{xA} [e^{2\pi A} v-v], \quad x\in[0,2\pi].
\lb{2.2}
\end{align}
Let us denote  $a=\max \{| \alpha' (x) | : x\in[0,2\pi]\}$ and
$b=\max \{\|e^{xA}\|: x\in [0,2\pi]\}$.
Note, that $\|e^{2\pi A}v\| =\| e^{(2\pi
-x)A} e^{xA} v\| \leq b\|e^{xA}v\|$ for any $x\in [0,2\pi]$.

First let us suppose that $F=L_p ([0,2\pi];E)$.
Then
$$\|Bg\|_{L_p ([0,2\pi];E)}
\leq (2\pi)^{1/p} \frac{ab}{m}.$$
On the other hand,
$$
\|g\|^p_{L_p ([0,2\pi];E)}
\geq \int_{\frac{4\pi}3}^{2\pi} \|e^{xA}v\|^p dx
\geq \dfrac{2\pi}3 b^{-p} \|e^{2\pi A} v\|^p \geq \dfrac{2\pi}{3} b^{-p}
(1-\dfrac1m)^p.
$$
Finally,
\begin{equation}
\|Bg\|_{L_p([0,2\pi];E)} \leq (2\pi)^{1/p}\dfrac{ab}{m} \leq
3^{1/p} ab^2 \|g\|_{L_p ([0,2\pi];E)} \cdot
\dfrac1{m-1}.
\lb{est}
\end{equation}
Since this holds for all $m$, it follows that $0\in \sig_{ap}(B)$.

Now suppose that $F=C([0,2\pi], E)$. Then
\begin{equation}
\|Bg\|_{C([0,2\pi],E)} \leq \dfrac{ab}{m}, \; \|g\|_{C([0,2\pi],E)} \geq
\|g(0)\|_E =\| e^{2\pi A}v\| \geq 1-\dfrac1m
\lb{bcontest}
\end{equation}
and hence
\begin{equation}
\|Bg\|_{C([0,2\pi];E)} \leq \dfrac{ab}{m-1} \|g\|_{C([0,2\pi];E)}.
\lb{contest}
\end{equation}
Since this is true for all $m$, it follows that
$0\in \sig_{ap} (B)$.
\end{pf}

\begin{thm}\lb{t2.1}
Let
$F$ be one of the spaces
$L_p ([0,2\pi],E)$, $1\leq p < \infty$
or  $C([0,2\pi],E)$.
Then the following are equivalent:
\newcounter{aa}
\begin{list}{\arabic{aa})}{\usecounter{aa}}
\item $1\in \rho (e^{2\pi A} )$ in $E$;
\item $1\in \rho (e^{2\pi B})$ in $F$;
\item $0\in\rho(B)$ in $F$.
\end{list}
\end{thm}

\begin{pf}
1) $\Rightarrow$ 2). Note that $(e^{2\pi B}f)(x) =e^{2\pi A} f(x)$.  Hence
$\sig (e^{2\pi A}; E)=\sig (e^{2\pi B}; F)$.
Note also that $\sig_r (e^{2\pi
A}; E)=\sig_r (e^{2\pi B}; F)$.

2) $\Rightarrow$ 3) follows from the spectral inclusion theorem
$e^{2\pi\sig (B)}\subset \sig(e^{2\pi B})$ (see \cite[p. 45]{24}).

3) $\Rightarrow$ 1).  Assume $0\in \rho (B; F)$ but
$1 \in \sig (e^{2\pi A}; E)
=\sig_{ap} (e^{2\pi A}; E) \cup \sig_r (e^{2\pi A}; E)$.
By Lemma~\ref{l2.1} it follows that
$1\in \sig _r (e^{2\pi A}; E)$, and hence
$1\in\sig_r (e^{2\pi B}; F)$.  By the spectral
mapping theorem for the residual spectrum (\cite[Theorem~2.5~(ii)]{24})
it follows that $ik\in\sig_r (B; F)$ for some $k\in\bbZ$.
Since $\sig(B; F)$ is
invariant under translations by $i$, we have that $0\in \sig (B; F)$,
contradicting 3).
\end{pf}

\section*{2.2. Spectral Mapping Theorem for Banach Spaces}\label{sub2.2}

As it is well-known (see, e.g., \cite[p.~82--89]{21}),
that in general, the inclusion
$e^{t \sig(A)} \subset \sig (e^{At}) $, $t\neq 0$
for a
$C_0$-semigroup $\{e^{tA}\}_{t\geq 0}$ on a Banach space $E$
is improper. In
particular, $i\bbZ \subset \rho (A)$ is implied by but does not imply
$1\in \rho (e^{2\pi A})$.  For Hilbert space $E$, however, the following
spectral mapping theorem of L.~Gearhart (see \cite[p.~95]{21}) is true:
{\it $1\in \rho (e^{2\pi A})$ if and only if $i\bbZ \subset \rho (A)$ and
$\sup_{k\in \bbZ} \| (A-ik)^{-1}\|< \infty$\/}.
We will now give a direct
generalization of this result to any arbitrary Banach space $E$.
This generalization
is related (but independent) to G.~Greiner's spectral mapping
theorem \cite[p.~94]{21} that involves C\'esaro summability of the
series $\sum_k (A-ik)^{-1}v$, $v\in E$.

\begin{thm} \lb{t2.2}
Let $\{e^{tA}\}$ be any $C_0$-semigroup on a Banach space $E$ and
let $F$ be one of the spaces
$L_p ([0,2\pi],E)$, $1\leq p < \infty$
or  $C([0,2\pi],E)$.
Then the following are equivalent:
\begin{list}{\arabic{aa})}{\usecounter{aa}}
\item $1\in \rho (e^{2\pi A})$;
\item $i\bbZ \subset \rho (A)$ and there is a constant $C>0$ such that
\begin{equation}
\|\sum_k (A-ik)^{-1}e^{ikx} v_k\|_F \leq C\| \sum_k
e^{ikx} v_k\|_F
\lb{G}
\end{equation}
for
any finite sequence $\{v_k\}\subset E$.
\end{list} \end{thm}

\begin{pf}
Consider the evolutionary semigroup
$\{e^{tB}\}_{t\geq 0}$ on $F$ from the previous subsection.  Consider a finite
sequence $\{v_k\}\subset E$. Assume that $(A-ik)^{-1}$ exists for all
$k\in \bbZ$. Define
functions $f$, $g\in F$ by the rule
\begin{equation}
f(x) =\sum_k (A-ik)^{-1} e^{ikx} v_k, \; g(x) =\sum_k e^{ikx} v_k,\quad
x\in [0,2\pi].
\lb{fg}
\end{equation}
Since $(A-ik)^{-1}: E\to \calD (A)$, one has $Bf=g$.  Indeed
\begin{align*}
(Bf)(x) &= \dfrac{d}{dt} e^{tA} f([x-t] (\mod 2\pi))\big|_{t=0}\\
& =\sum_k [A (A-ik)^{-1} e^{ikx} v_k -ik (A-ik)^{-1} e^{ikx}v_k]=g.
\end{align*}

1) $\Rightarrow$ 2).  If $1\in \rho (e^{2\pi A})$, then the
inclusion $i\bbZ \subset\rho (A)$
follows from the spectral inclusion theorem $e^{2\pi \sig(A)}\subset
\sig (e^{2\pi A})$.
In accordance
with part 1) $\Rightarrow$ 3) of Theorem~\ref{t2.1},
the operator $B$ has bounded
inverse $B^{-1}$ on $F$ provided that
$1\in \rho (e^{2\pi A})$.  Denote $C=\|B^{-1}\|$, and consider
functions \eqref{fg}.  Then $\|f\|_F=\|B^{-1}g\|_F \leq C\|g\|_F$,
and \eqref{G} is proved.

2) $\Rightarrow$ 1).  First, we show that 2) implies $0\notin \sig_{ap} (B)$.
Indeed, the
functions of type $g$ in \eqref{fg} are dense in $F$.
If we let $u_k =(A-ik)^{-1} v_k$, then we note that the functions
of type $f$ are also dense in $F$. Now \eqref{G} implies
$\|Bf\|_F =\|g\|_F \geq C^{-1} \|f\|_F$, and $0 \notin
\sig_{ap}(B)$.

Assume that 2) is fulfilled, but $1\in \sig (e^{2\pi A})=\sig_r (e^{2\pi A})
\cup \sig_{ap} (e^{2\pi A})$.  If $1\in \sig_{ap}(e^{2\pi A})$ in $E$
then, by Lemma~\ref{l2.1}, $0\in \sig_{ap}(B)$, in contradiction to the
previous paragraph.  On the other hand,
$1\in \sig_r (e^{2\pi A})$ implies, by the
spectral mapping theorem for residual spectrum, that $ik \in \sig_r (A)$
for some $k\in\bbZ$, contradicting
$i\bbZ \subset \rho (A)$.
\end{pf}
\noindent{\bf Remark.} We note, that 1) $\Rightarrow$ 2) can be also
seen directly. Indeed, assuming 1), let us denote
$\phi (s)=(e^{2\pi
A}-I)^{-1}e^{sA}$, $s\in [0,2\pi]$. Then the convolution operator
\begin{equation*} (Kf)(x)=\int\limits_0^{2\pi}\phi(s)f(x-s)\,ds
\end{equation*}
is a bounded operator on $F$.
But
\begin{equation*}
(A-ik)^{-1}=\int\limits_0^{2\pi}e^{-iks}(e^{2\pi A}-I)^{-1}e^{sA}\,
ds,\quad k\in\Bbb Z,
\end{equation*}
are Fourier coefficients of $\phi:[0,2\pi]\to B(E)$. Inequality
\eqref{G} can be viewed as the condition of boundedness of $K$,
which gives 2).
\medskip

Let us show now that Theorem~\ref{t2.2} is really a direct generalization
of L.~Gearhart's Theorem, mentioned above.  Indeed, for Hilbert space $E$ and
$p=2$, Parseval's identity implies:
\begin{align*}
& \|\sum_k (A-ik)^{-1} e^{ik\cdot} v_k\|_{L_2 ([0,2\pi];E)} =\Big( 2\pi
\sum_k \|(A-ik)^{-1} v_k \|_E^2\Big)^{1/2}\\
& \| \sum_k e^{ik\cdot} v_k \|_{L_2 ([0,2\pi];E)} =\Big( 2\pi \sum_k
\|v_k\|_E^2\Big)^{1/2}.
\end{align*}
Clearly, \eqref{G} is equivalent to the condition
$\sup \{\|(A-ik)^{-1} \|:
k\in\bbZ \} < \infty$.

We conclude this subsection by giving four more statements equivalent
to 1) and 2) in Theorem~\ref{t2.2}:

\begin{list}{\arabic{aa})}{\usecounter{aa}}\setcounter{aa}{2}
\item $i\bbZ \subset \rho (A)$ and there is a constant $C>0$ such that
$$
 \|Bf\|_{L_1([0,2\pi]; E)}\geq
 C^{-1}\|f\|_{C([0, 2\pi];E)}$$ for all  $f\in
C([0, 2\pi]; E)$ such that $ Bf\in L_1([0,2\pi]; E)$;
\item $i\bbZ \subset \rho (A)$ and there is a constant $C>0$ such that
$$\|Bf\|_{C([0, 2\pi];E)}\geq
 C^{-1}\|f\|_{L_1([0,2\pi]; E)}$$ for all
$f\in L_1([0,2\pi]; E)$ such that  $Bf\in C([0, 2\pi]; E);$
\item $i\bbZ \subset \rho (A)$ and there is a constant $C>0$ such that
\begin{equation*}
\|\sum_k (A-ik)^{-1}e^{ikx} v_k\|_{C([0, 2\pi];E)}
\leq C\| \sum_k
e^{ikx} v\|_{L_1([0,2\pi]; E)}
\end{equation*}
for any finite sequence $\{v_k\}\subset E$;
\item $i\bbZ \subset \rho (A)$ and there is a constant $C>0$ such that
\begin{equation*}
\|\sum_k (A-ik)^{-1}e^{ikx} v_k\|_{L_1([0,2\pi]; E)}
\leq C\| \sum_k
e^{ikx} v\|_{C([0, 2\pi];E)}
\end{equation*}
for any finite sequence $\{v_k\}\subset E$.
\end{list}

\section*{2.3.  Real Line} \label{sub2.3}
Consider now the evolutionary
semigroup $\{e^{tB}\}_{t\geq 0}$,
\begin{equation}\label{evolaut}
(e^{tB} f)(x)
=e^{tA}f(x-t)
\end{equation}
acting on the space $F=L_p (\bbR;E)$, $1\leq p <\infty$,
or $F=C_0 (\bbR;E)$.  Formula
\eqref{genB} is valid and the identities
\begin{align*}
& e^{tB}e^{i\xi\cdot } f(\cdot )=
e^{i\xi\cdot } e^{-i\xi t} e^{tB} f(\cdot ), \\
& Be^{i\xi\cdot }f(\cdot) =e^{i\xi \cdot} (B- i\xi ) f(\cdot),\quad \xi
\in\bbR,
\end{align*}
show that
$\sigma (e^{tB})$ in $F$ is invariant under rotations
centered at the origin, and that
$\sig(B)$, $\sig_{ap} (B)$ and $\sig_r (B)$ in $F$
are invariant under translations parallel to $i\bbR$.

First we state a simple lemma. Let $F_s$ be one of the spaces
$F_s=l_p(\bbZ;E)$, $1\leq p<\infty$ or $F_s=c_0(\bbZ;E)$ (of
sequences $(v_n)_{n\in \bbZ}$ such
that  $v_n\to 0$ as $n\to\pm\infty$). Let $S$ be
a shift operator on $F_s$, that is $S: (v_n)\mapsto (v_{n-1})$.
For an operator $a$ on $E$ we will denote by $D_a$  the
diagonal operator
on $F_s$, acting by the rule
$D_a:(v_n)_{n\in\bbZ} \mapsto (av_n)_{n\in\bbZ}$.

\begin{lem}\label{aS}
The following are equivalent:
\begin{list}{\arabic{aa})}{\usecounter{aa}}
\item $\sig (a) \cap \bbT =\emptyset$ in $E$;
\item $\sig (D_aS) \cap \bbT =\emptyset$ in $F_s$.
\end{list}
\end{lem}

\begin{pf}
We will give the proof for the case when $F_s = l_p=l_p(\bbZ;E)$.  The case
$F_s=c_0(\bbZ;E)$ can be considered similarly.

1) $\Rightarrow$ 2).  Since $\sig (a) \cap
\bbT =\emptyset$, there exists a Riesz projection ${\hat p}$ for
$a$ in $E$ that corresponds to the part of the spectrum $\sig
(a) \cap\bbD$. Define ${\hat q}=I-{\hat p}$, and consider in $F$
complimentary projections $D_{\hat p}$ and $D_{\hat q}$. Since
$D_aSD_{\hat p}=D_aD_{\hat p}S$, the decomposition
$F_s=\Im D_{\hat p} \oplus \Im D_{\hat q}$ is
$D_aS$-invariant.
For the
spectral radius
$$r(\cdot)= \lim_{n\to \infty} \|(\cdot )^n\|^{\frac{1}{n}}$$
one has $r({\hat p}a{\hat p})<1$ and $r([{\hat q}a{\hat q}]^{-1})<1$
in $E$. Hence,
$r(D_aSD_{\hat p})<1$, $r([D_aSD_{\hat q}]^{-1})<1$,
and $\sig (D_aS) \cap \bbT =\emptyset$ in $F_s$.

2) $\Rightarrow$ 1).
Assume that $I-D_aS$ is invertible in $l_p$, but
for any $\epsilon >0$ there is a vector $v\in E$ such that
$\|v\|_E=1$ and $\|v-av\|_E < \epsilon$. Fix $q>0$ such that
$$
\Big| 1-e^{\pm q}\Big|  < \epsilon,
$$
and define a sequence $(v_n)\in l_p$ by $v_n=e^{-q|n|}v$,
$n\in \bbZ$. Then
$$
I-D_aS : (v_n)_{n\in\bbZ}
\mapsto \Big( e^{-q|n|}(v-av) +
(e^{-q|n|}- e^{-q|n-1|})av \Big)_{n\in\bbZ}.
$$
A direct calculation shows that
$$
\|(I-D_aS)(v_n)\|_{l_p}\leq (1+\|a\|) \cdot \epsilon \cdot \|(v_n)\|_{l_p},
$$
contradicting the invertibility of $I-D_aS$ in $l_p$.

Let us show now that $I-a$ has a dense range in $E$,
provided $I-D_aS$ is an operator onto $l_p$. Indeed, for any
$u\in E$ consider a sequence $(u_n)\in l_p$ defined by $u_0=u$ and
$u_n=0$ for $n\neq 0$. Find a sequence $(v_n)\in l_p$ such that
$(I-D_aS)(v_n)=(u_n)$, that is $v_n-av_{n-1}=u_n$ for $n\in \bbZ$.
But then for  $k\in \bbN$ one has
\begin{align*}
& u= \sum_{n=-k}^k(v_n-av_{n-1}) \\
& = v_k -av_{-k-1} + (y_k - ay_k) ,
\end{align*}
where
$$ y_k= \sum_{n=-k}^{k-1}v_n.$$
Therefore, $\Im (I-a) \ni y_k-ay_k\to u$, since
$v_k\to 0$  and $av_{-k-1}\to 0$ as $k\to\infty$.
\end{pf}

\begin{thm}\lb{t2.3}
Let $F$ be one of the spaces
$L_p (\bbR;E)$, $1\leq p < \infty$
or  $C_0(\bbR;E)$, and let $t>0$.
Then the following are equivalent:
\begin{list}{\arabic{aa})}{\usecounter{aa}}
\item $\sig (e^{tA}) \cap \bbT =\emptyset$ in $E$;
\item $\sig (e^{tB}) \cap \bbT =\emptyset$ in $F$;
\item $0\in \rho (B)$ in $F$.
\end{list} \end{thm}

\begin{pf}
2) $\Rightarrow$ 3) follows from the spectral inclusion theorem
for $\{e^{tB}\}$.

3) $\Rightarrow$  2)
we will prove for $F=L_p(\bbR;E)$;
the arguments for $F=C_0(\bbR;E)$ are similar.

Since $\sigma (e^{tB})$ is invariant
under the rotations with the center at origin, it
suffices to prove that 3) implies
$1\in \rho (e^{tB})$. Also, to confirm our
previous notations, we will consider only the case $t=2\pi$.
The proof stays the same for any $t$.

The idea is to apply Theorem~\ref{t2.1}, to show that
$1\in\rho(e^{2\pi B})$ is implied by
$0\in\rho(B')$. Here the operator
$B'= -{d\over ds} - {d\over dx} + A$ acts on
$L_p([0,2\pi]\times \bbR; E)$, $s\in[0,2\pi]$,
$x\in \bbR$.
Indeed, by formula \eqref{genB} one has $B=-\frac{d}{dx}+A$.
Hence, $B'$ on $L_p([0,2\pi]; L_p(\bbR; E))$ is the generator
of the evolutionary semigroup for the semigroup
$\{e^{tB}\}$ on $L_p(\bbR; E)$.
But the change of variables $u = [s-x] (\mod 2\pi)$,
$v = x$ shows that $\rho(B') = \rho (-{d\over dv} + A) = \rho(B)$.
Let us now make this argument more formal.

Consider the semigroups
\begin{align*}
( e^{tB'}h)(s,x) &= e^{tA} h( [s-t] (\mod 2\pi),\, x-t),~~t>0,
~s\in [0,2\pi],~x\in \bbR, \\
(e^{t\calB}h)(s,x) &= e^{tA}h(s,\,x-t),~~t>0,
~s\in [0,2\pi],~x\in \bbR,
\end{align*}
and an invertible isometry $J$,
$$
(Jh)(s,x)= h( [s+x] (\mod 2\pi),\, x),
$$
acting on the space
$$L_p([0,2\pi]\times \bbR;\, E) = L_p([0, 2\pi]; L_p(\bbR; E)).
$$
Since $e^{t\calB}$ acting
on $L_p([0, 2\pi]; L_p(\bbR; E))$ is actually the
operator of multiplication by the operator $e^{tB}$ in $L_p(\bbR;E)$,
one has:
$$\sigma (e^{t\calB})=\sigma (e^{tB})\mbox{ and }
\sigma(\calB; L_p([0,2\pi]; L_p(\bbR; E)))=\sigma(B; L_p(\bbR; E)).$$
Also,
\begin{align*}
& \Big(Je^{tB'}h\Big)(s,x)=
e^{tA}h( [s+x-t](\mod 2\pi), x-t )=  \\
& \Big(e^{t\calB}Jh\big)(s,x),
\end{align*}
and hence
one has $Je^{tB'}=e^{t\calB}J$ and $JB'=\calB J$. Therefore,
$$\sigma (e^{t\calB})=\sigma (e^{tB'}) \mbox{ and }
\sigma(\calB)=\sigma(B')\mbox{ in } L_p([0,2\pi]; L_p(\bbR; E)).$$
Thus 3) implies
$0\in\rho (\calB)$ and $0\in \rho (B')$.

The semigroup $\{e^{tB'}\}$ acts on $L_p([0, 2\pi]; L_p(\bbR; E))$
by the rule
$$(e^{tB'}f)(s)=e^{tB}f( [s-t](\mod 2\pi)),$$
where $f(s)=h(s, \cdot)\in L_p(\bbR; E)$ for almost all
$s\in [0,2\pi]$.
Hence, $\{e^{tB'}\}$
on the space $L_p([0,2\pi]; L_p(\bbR; E))$
is the evolutionary semigroup for the semigroup $\{e^{tB}\}$
on $L_p(\bbR; E)$. Now
one can apply the part 3) $\Rightarrow$ 1) of the Theorem~\ref{t2.1}
and conclude that $1\in\rho (e^{2\pi B})$ on $L_p(\bbR;E)$.

1) $\Leftrightarrow$ 2)
we will prove for $F=L_p(\bbR;E)$;
the arguments for $F=C_0(\bbR;E)$ are similar.

Let us denote, for brevity,
$a=e^{2\pi {\cal A}}$ and $(Rf)(x)=f(x-2\pi)$ on $L_p(\bbR;E)$.
Thus $e^{2\pi B}=aR$. Consider the invertible isometry
\begin{align*}
& j: L_p(\bbR;E)\to l_p(\bbZ; L_p([0, 2\pi];E)) : \\
& f\mapsto (f_n),~~f_n(s)=f(s+2\pi n),~n\in\bbZ, s\in [0,2\pi).
\end{align*}
Let $S: (f_n)\mapsto (f_{n-1})$ be a shift
operator on  $l_p(\bbZ; L_p([0, 2\pi];E))$.
Then $jaR=D_aSj$ and  $\sigma (aR) = \sigma (D_aS)$. Therefore,
2) is equivalent to $\sig (D_aS) \cap \bbT =\emptyset$.
By Lemma~\ref{aS}
this in turn is equivalent to  $\sig (a) \cap \bbT =\emptyset$ in
$L_p(\bbR; E)$.
\end{pf}

Note, that 3) $\Rightarrow$ 1) in the above theorem
can also be derived directly by constructing a function $g$ in a similar
manner as in the proof of Lemma~\ref{l2.1}. This proof will be given
elsewhere.

\section{Non-autonomous Case}

Consider a non-autonomous differential equation $v'(x) =A(x) v(x)$,
$x\in\bbR$ in $E$.  We will assume that this equation is well-posed.
This means that there exists an evolutionary family $\{U(x,s)\}_{x\geq s}$
(propagator) for the equation, that is $v(x) =U(x,s) v(s)$, $x\geq s$.
Recall the definition of evolutionary family (see, e.g., \cite{27,31}).
\begin{defn}\label{evfam}
A family $\{U(x,s)\}_{x\geq s}$ of bounded in $E$ operators
$U(x,s)$ is called an {\it evolutionary family} if the following conditions
are
fulfilled:\\
(i)  for each $v\in E$ the
function $(x,s) \mapsto U(x,s) v$ is continuous for
$x\geq s$;\\
(ii) $U(x,s) =U(x,r) U(r,s)$, $U(x,x) =I$, $x\geq r \geq s$;\\
(iii) $\|U(x,s)\|\leq Ce^{\beta (x-s)}$, $x\geq
s$ for some constants $C$, $\beta >0$.
\end{defn}

The evolutionary family $\{U(x,s)\}$ generates an
evolutionary semigroup $\{e^{tG}\}_{t\geq 0}$ acting on the space
$F=L_p(\bbR;E)$, $1\leq p <\infty$
or $F=C_0 (\bbR;E)$
by the rule
\begin{equation}\label{evnon}
(e^{tG}f)(x) =U(x,x-t) f(x-t),~~x\in\bbR.
\end{equation}
In Subsection~3.1
we will prove the spectral mapping theorem $\sig(e^{tG})\bs \{0\}
=e^{t\sig (G)}$, $t\neq 0$ for $\{e^{tG}\}$.  We will
achieve this
by applying a simple change of variables argument
(cf.\ the proof of Theorem~\ref{t2.3}) to deduce this
from
Theorem~\ref{t2.3}.
In Subsection~3.2 we will prove
that the hyperbolicity $\sig (e^{tG}) \cap \bbT =\emp$ of the semigroup
in $F$ is equivalent to the so-called spectral hyperbolicity of the family
$\{U(x,s)\}$.  Spectral hyperbolicity is a generalization of the notion of
exponential dichotomy
(see, e.g.,
\cite{9})
for the equation $v'(x) =A(x) v(x)$ with bounded $A:\bbR\to B(E)$.

\section*{3.1.  The Spectral Mapping Theorem for Evolutionary Semigroup}

Let $G$ be the generator of the evolutionary semigroup $\{e^{tG}\}_{t\geq
0}$ acting on the space
$F=L_p(\bbR;E)$, $1\leq p <\infty$
or $F=C_0 (\bbR;E)$ by equation \eqref{evnon}.

\begin{thm}\lb{t3.1}
The spectrum $\sig (G)$ is invariant under translations along the
imaginary axis, and the following are equivalent:\\
1) $0\in \rho (G)$ on $F$;\\
2) $\sig (e^{tG}) \cap \bbT =\emp$ on $F$,
$t>0$.
\end{thm}

\begin{pf}
For any $\xi \in \bbR$ it is
true that $e^{tG} e^{i\xi\cdot} f(\cdot) =e^{i\xi
(\cdot -t)}e^{tG} f(\cdot)$ and $Ge^{i\xi \cdot} =e^{i\xi \cdot}
(G-i\xi)$.  Hence $\sig (e^{tG})$ is invariant under rotations
centered at the origin, and $\sig (G)$ is invariant under translations along
the imaginary axis.

2) $\Rightarrow$ 1) follows from the spectral inclusion theorem
for $\{e^{tE}\}$.

1) $\Rightarrow$ 2). We will first consider the case when
$F=L_p(\bbR;E)$, $1\leq p <\infty$.

The idea of the proof is almost identical to the
proof of 3) $\Rightarrow$ 2) from Theorem~\ref{t2.3}.
If $U(\cdot, \cdot)$ is a smooth propagator for the equation
$v'=A(x)v$, then $G=-\frac{d}{dx}+A(x)$. Consider
the evolutionary semigroup for
$\{e^{tG}\}$,
that is the semigroup with the generator
$B=-\frac{d}{ds}+G$ on $L_p(\bbR; L_p(\bbR; E))$.
Theorem~\ref{t2.3} shows that
$1\in\rho(e^{tG})$ is implied by
$0\in\rho(B)$, where
$B = -{d\over ds} - {d\over dx} + A(x)$, $s,x\in \bbR$
by formula \eqref{genB}.
The change of variables $u = s-x$,
$v = x$ shows that $\rho(B) = \rho (-{d\over dv} + A(v))
= \rho(G)$.  Let us now
make this argument more formal.

Consider the semigroups $\{e^{tB}\}_{t\geq 0}$ and $\{e^{t\calG}\}_{t\geq 0}$
acting on the space $L_p (\bbR \times \bbR; E) =L_p (\bbR;
L_p (\bbR;E))$ by
\begin{align*}
(e^{tB}h) (s,x) &=U(x,x-t) h(s-t, x-t), ~~(s,x) \in\bbR^2, ~t>0,\\
(e^{t\calG} h)(s,x) &=U(x,x-t) h(s,x-t).
\end{align*}
Note that $e^{t\calG}$ in $L_p(\bbR; L_p(\bbR; E))$
is the operator of multiplication by $e^{tG}$, that is
$(e^{t\calG}f)(s)=e^{tG}f(s)$, where $f(s)=h(s, \cdot )\in L_p(\bbR; E)$.
Similarly, $(\calG f)(s)=Gf(s)=Gh(s, \cdot )$ for
$f(s) = h(s, \cdot )\in \calD_{L_p(\bbR; E)}(G)$
for almost all $s\in\bbR$.

Consider an isometry $J$ on $L_p (\bbR\times \bbR;E)$ defined by
$(Jh)(s,x)=h(s+x, x)$.  Then for $h\in L_p (\bbR\times \bbR;E)$
one has:
\begin{equation}
(e^{t\calG} Jh) (s,x) =U(x,x-t) h(s+x-t, x-t)=(Je^{tB}h)(s,x).
\lb{cb}
\end{equation}
Also \eqref{cb} implies
$$
\calG Jh =JBh, \; h\in\calD(B)
\mbox{ and } J^{-1} \calG h =BJ^{-1} h, \;
h\in\calD (\calG).
$$
Therefore, $\sig (\calG) =\sig (B)$ on $L_p(\bbR \times \bbR; E)$.

Note that $\calG$ on $L_p(\bbR \times \bbR;E)$ has bounded inverse
$(\calG^{-1} f)(s)=G^{-1} f(s)$, $s\in\bbR$, provided $G$ has
bounded inverse $G^{-1}$ on $L_p (\bbR;E)$.  Hence 1) implies $0\in\rho
(B)$.

Let us apply the part  3) $\Rightarrow$ 1) of
Theorem~\ref{t2.3} to the semigroups $\{e^{tG}\}$ and
$\{e^{tB}\}$.  To this end we note that $(e^{tB}f)(s) =e^{tG} f(s-t)$ for
$f:\bbR \to L_p (\bbR;E): s\mapsto h(s,\cdot)$.  Hence $0\in\rho (B)$ on
$L_p (\bbR;L_p (\bbR;E))$ implies $\sig (e^{tG}) \cap \bbT =\emp$,
$t\neq 0$ on $L_p (\bbR;E)$.

The proof for the case $F= C_0(\bbR;E)$ is
identical, and uses exactly the same semigroups and isometries
on $C_0(\bbR;F) = C_0(\bbR\times\bbR;E)$.
\end{pf}

\section*{3.2.  Hyperbolicity}

Let $\{U(x,s)\}_{x\geq s}$ be an
evolutionary family on a Banach space $E$. In this subsection
we relate the (spectral) hyperbolicity of the evolutionary
family and the hyperbolicity of the evolutionary
semigroup $\{e^{tG}\}$ on the space $L_p=L_p(\bbR;E)$ in the case
when the
Banach space $E$ is separable.
The case $F= C_0(\bbR;E)$
(without the assumption of separability)
and the case of a Hilbert space $E$ and $p=2$
was considered in \cite{27,28}.

\begin{defn}\lb{hyp}
An evolutionary family $\{U(x,s)\}_{x \geq s}$ is called (spectrally)
hyperbolic if there exists a projection-valued, bounded function $P:\bbR
\to B(E)$ such that the function $\bbR \ni x \mapsto P(x) v\in E$ is
continuous for every $v\in E$ and for some constants $M$, $\lam >0$
and all $x\geq s$ the
following conditions are fulfilled:
\newcounter{cc}
\begin{list}{\alph{cc})}{\usecounter{cc}}
\item  $P(x) U (x,s) =U(x,s) P(s)$:
\item  $\| U(x,s) v\| \leq Me^{-\lam (x-s)} \|v\|$ if $v\in \Im P(s)$,\\
$\|U(x,s) v \| \geq M^{-1} e^{\lam (x-s)} \|v\|$ if $v\in \Ker P(s)$;
\item $\Im (U(x,s) | \Ker P(s))$ is dense in $\Ker P(x)$.
\end{list}
\end{defn}

This notion generalizes the notion of exponential dichotomy
(see, e.g.,
\cite{9}) for the
solutions of differential equation $v'(x) =A(x) v(x)$, $x\in\bbR$,
with bounded and continuous $A:\bbR \to B(E)$.
In this case the evolutionary family (propagator)
$\{U(x,s)\}_{(x,s)\in\bbR^2}$ consists of invertible operators,
the function $(x,s)
\mapsto U(x,s)$ is norm-continuous, and $P(\cdot)$ from
Definition~\ref{hyp} is also a
bounded, norm-continuous function $P:\bbR \to
B(E)$.

The second inequality in b) implies that
the restriction $U(x,s) | \Ker P(s)$ is
uniformly injective
as an operator from $\Ker P(s)$ to $\Ker P(x)$
(that is $\|U(x,s) v\| \geq c\|v\|$
for some $c> 0$
and all $v\in\Ker P(s)$).
Thus condition~c)
implies that $U(x,s) | \Ker P(s)$ is invertible as an operator from
$\Ker P(s)$ to $\Ker P(x)$.  Obviously, if $U(x,s)$ is
invertible in $E$ or
$\dim \Ker P(x) \leq d < \infty$, condition~c) in Definition~\ref{hyp} is
redundant. The inequality
$\dim \Ker P(x) \leq d < \infty$ holds,
for example, if the $U(x,s)$ are compact
operators in $E$ ([R.~Rau, private communication]).
Generally, of course, b) does not imply c).

>From now on we will assume that the Banach space $E$ is separable.

As we will see below,
the spectral hyperbolicity of the evolutionary
family $\{U(x,s)\}$ is equivalent
to the hyperbolicity $\sig
(e^{tG} ) \cap \bbT =\emp$, $t> 0$ of the evolutionary semigroup
$\{e^{tG}\}_{t\geq 0}$ in $L_p (\bbR;E)$.
Therefore, by Theorem~\ref{t3.1}
the spectral hyperbolicity of the evolutionary
family $\{U(x,s)\}_{x\geq s}$ is also equivalent
to the condition $\sig (G) \cap i\bbR =\emp$.
That is why we used the term
{\it spectral} hyperbolicity in Definition~\ref{hyp}.  A remarkable
observation by R.~Rau \cite{27} shows that generally
the condition c) in
Definition~\ref{hyp} cannot be dropped: there exists an evolutionary
family that satisfies conditions a) and b) but $\sig (e^{tG})\cap \bbT
\neq \emp$ for the associated evolutionary semigroup.

If the operator
$T=e^G$ is hyperbolic in $L_p (\bbR;E)$, that is $\sig (T) \cap \bbT
=\emp$, we let $\calP$ denote the Riesz projection for $T$, corresponding to
the part $\sig (T)$ lying
inside the unit disk $\bbD$, and set $\calQ =I-\calP$.

\begin{lem}\lb{ling}
If $\sig (e^{tG})\cap \bbT =\emp$, $t> 0$, then $\calP$ has a form
$(\calP f)(x) =P(x)f(x)$, where $P:\bbR \to B(E)$ is a bounded
projection-valued function such that
the function $\bbR \ni x \mapsto P(x) v \in E$ is
(strongly) measurable for each $v\in E$.
\end{lem}

\begin{pf}
We will show first that
\begin{equation}
\chi \calP =\calP\chi
\lb{com}
\end{equation}
for any scalar function $\chi \in L_\infty (\bbR; \bbR)$.
Then we will derive
that $(\calP f)(x) =P(x)f(x)$ from \eqref{com}.

Note that the
decomposition $L_p(\bbR;E)=\Im \calP \oplus \Im \calQ$ is $T$-invariant.
Denote $T_P =\calP T\calP =T| \Im \calP$,
$T_Q =\calQ T\calQ=T| \Im \calQ$.  Note that
$\sig (T_P) \subset \bbD$, and $T_Q$ is invertible with $\sig
(T^{-1}_Q) \subset \bbD$ in $\Im \calQ$.
Hence for some $\lam$, $M> 0$ and all $n\in \bbN$,
the following inequalities hold:
\begin{align}
\| T_P^n f\|_{L_p} & \leq Me^{-\lam n} \|f\|_{L_p}, \quad f\in \Im \calP,
\lb{est1}\\
\| T_Q^n f\|_{L_p} & \geq M^{-1} e^{\lam n} \| f\|_{L_p}, \quad f\in \Im
\calQ.
\lb{est2}
\end{align}

We show first that $\Im \calP =\{f\in L_p (\bbR;E): T^n f\to 0$ as $n\to
\infty\}$.  Indeed, $f\in \Im \calP$ implies that
$T^n f\to 0$ by \eqref{est1}.
Conversely, if $T^n f\to 0$, then for
$f=\calP f +\calQ f$, the inequality \eqref{est2} implies
$$
\| \calQ f\| \leq Me^{-\lam n} \| T_Q^n \calQ f\|
\leq Me^{-\lam n} \{\|T^n f\| +\|
T^n_P f\|\}\to 0,
$$
and hence $f\in \Im \calP$.

Consider on $L_p=L_p(\bbR;E)$ the
operator $\chi$ of multiplication by $\chi (\cdot ) \in L_\infty
(\bbR;\bbR)$.  Note that $(T^n \chi f)(x) =\chi (x-n)(T^n f)(x)$. Hence
for $f\in \Im \calP$
$$
\| T^n \chi f\|_{L_p}
\leq \| \chi\|_{L_\infty} \| T^n f\|_{L_p} \to 0 \mbox{ as } n
\to \infty,
$$
and so $\chi f \in \Im \calP$.  Thus,
to prove \eqref{com}, it suffices to show that
$f\in \Im \calQ$ implies $\chi f \in \Im \calQ$.

Fix $f\in \Im \calQ$.  Recall that $T_Q$ is invertible on $\Im \calQ$.  Let
$f_n =T^{-n}_Q f$, and define functions $g_n (x) =\chi (x+n) f_n
(x)$, $n=0,1,\ldots$ .  Decompose $g_n =\calP g_n +\calQ g_n$.
Since the decomposition
$L_p(\bbR;E) =\Im \calP \oplus \Im \calQ$ is $T$-invariant, one has:
$$
\chi f =T^n g_n, \; \calP\chi f =T^n_P \calP g_n, \;
\calQ\chi f =T^n_Q \calQ g_n, \;
n=0,1,\ldots.
$$
Now \eqref{est1}--\eqref{est2} imply:
\begin{align*}
&\| \calP\chi f\| \leq Me^{-\lam n} \| \calP g_n\|
\leq Me^{-\lam n} \{ \| g_n\|
+\| \calQ g_n\| \} \\
& \leq Me^{-\lam n} \{ \| \chi \|_{L_\infty} \|f_n \| +Me^{-\lam n}
\|\calQ\chi f\|\}\\
& \leq Me^{-\lam n} \{ \| \chi\|_{L_\infty} Me^{-\lam n} \|f\|+Me^{-\lam n}
\|\calQ \chi f\| \} \to 0,
\end{align*}
and hence $\chi f \in \Im \calQ$. Thus \eqref{com} is proved.

In order to define $P(\cdot)$ such that
$(\calP f)(x) =P(x) f(x)$, fix $m\in \bbZ$ and let $\chi_m
(x)=1$ if $x\in [m,m+1)$, and $\chi_m (x) =0$ otherwise.
Let $\{e_n\}_{n\in\bbZ}$ be a linearly independent set with dense
span $E_0$.

Consider the function
$f\in L_p (\bbR;E)$, defined by $f(x) =\chi_m (x) e_n$.
Since $\calP$ is
bounded on $L_p(\bbR;E)$, it is true that $\calP f\in L_p (\bbR;E)$.
For
$x\in [m, m+1)$ define
a vector $P(x) e_n\in E$ as $P(x) e_n =(\calP f)(x)$.
For $v=\sum_{n=1}^k d_n e_n \in E_0$
set $P(x) v=\sum_{n=1}^k d_n P(x)
e_n$.

Let $\Del$ be a measurable subset in $[m, m+1)$, and let $\chi_\Del$ be its
characteristic function.  Now  \eqref{com} implies that
\begin{align*}
\int_\Del \| P(x) v\|_E^p \, dx &=
\int_{\bbR} \| \chi_\Del (\calP\chi_m v)(x)
\|_E^p \, dx =
\int_{\bbR} \| (\calP\chi_\Del v)(x) \|_E^p \, dx \\
& \leq \| \calP \|_{B(L_p (\bbR;E))}^p \int_\Del \|v\|_E^p \, dx.
\end{align*}
Therefore, $\|P(x) v\|_E \leq \|\calP\|_{B(L_p)} \|v\|_E$
for a.e.\ $x\in \bbR$ and all
$v\in E_0$.  Hence, $P(x)$ can be extended to a bounded operator on $E$,
such that $\|P(x) \| \leq \| \calP\|$ for a.e. $x\in\bbR$.  That
the function $x\mapsto P(x) v$ is
a measurable function
for all $v\in E_0$ (and, hence, for all $v\in E$) follows from the fact that
the function $x\mapsto
(\calP f) (x)$ is measurable.

To show that $(\calP f)(x) = P(x)f(x)$ we can assume that $f$ is
a simple function, $f=\sum\chi_{\Delta_k}v_k$, where
$\Delta_k \subset [m_k, m_k+1)$, $m_k\in \bbZ$,
and $v_k\in E_0$. Then \eqref{com}
implies:
$$(\calP f)(x)=\sum \chi_{m_k}(x)(\calP \chi_{\Delta_k}v_k)(x)=
\sum\chi_{\Delta_k}(x)P(x)v_k = P(x)f(x).
$$
\end{pf}

Let us stress that the function $P(\cdot )$ above is only defined
on a set $\bbR_0\subset\bbR$ such that $\mes (\bbR\setminus\bbR_0) = 0$.
In Theorem~\ref{cont} below we will establish that, in fact,
this function
$P(\cdot)$ can be extended to all of $\bbR$ as
a continuous function (in the strong
operator topology in $B(E)$).  To prove
this fact we will need a few definitions and Lemma.

Let $\frakA$ be the set of all operators $a$ in $L_p(\bbR;E)$ of
the form $(af)(x) =a(x) f(x)$,
where the function $a:\bbR \to B(E)$ is bounded
and the function
$\bbR \ni x \mapsto a(x) v\in E$ is continuous for each $v\in E$.
For $a\in\frakA$ and $x\in\bbR$ let us define an
operator $\pi_x (a)$ on $l_p (\bbZ;E)$ by the rule
\begin{equation}
\pi_x (a) =\diag \{a(x+n)\}_{n\in\bbZ}: (v_n)_{n\in\bbZ} \mapsto (a(x+n)
v_n)_{n\in\bbZ}.
\lb{pi}
\end{equation}
Finally, let $S: (v_n)_{n\in\bbZ} \mapsto (v_{n-1})_{n\in\bbZ}$ be a
shift operator on $l_p (\bbZ;E)$.

Let us denote: $T=e^G$, $a(x) =U(x,x-1)$, $(Rf) (x) =f(x-1)$.
Then $T=aR$. For $\lam \in\bbT$ set $b=\lam I -aR$, and for
$x\in\bbR$ set $\pi_x (b) =\lam I-\pi_x(a)S$.

\begin{lem}\lb{inv}
If $\sig (T) \cap \bbT =\emp$ in $L_p(\bbR;E)$ then $\sig (\pi_x (a) S)
\cap \bbT =\emp$ in $l_p(\bbZ;E)$
for all $x\in\bbR$.
Moreover, for all $\lam \in\bbT$ the following estimate holds:
\begin{equation}
\| [\pi_x (b)]^{-1}\|_{B(l_p (\bbZ;E))} \leq \|b^{-1}\|_{B(L_p
(\bbR;E))},\quad x\in\bbR.
\lb{binv}
\end{equation}
\end{lem}

\begin{pf}
First, for any $\xi \in\bbR$ one has:
$$\pi_x (a) SL=e^{-i\xi} L\pi_x (a) S,$$
where $L$ is the operator $(v_n)\mapsto (e^{i\xi n} v_n)$.
Hence
$\sig (\pi_x(a) S)$ is invariant under rotations centered at
the origin.  Thus it suffices to prove the Lemma for the
special case $\lam=1$, that is,
to show that if $b=I-aR$ is invertible in $L_p (\bbR;E)$
then for
each $x_0\in\bbR$ that the operator
$\pi_{x_0}(b)=I-\pi_{x_0}(a)S$ is invertible in $l_p(\bbZ;E)$, and that
estimate \eqref{B}
is valid for this $b$ and $x=x_0$.

Further, it suffices to prove the Lemma only for $x_0=0$.
Indeed, let us denote $\hat{a}(x)=a(x+x_0)$, $x\in \bbR$ for any fixed
$x_0\in \bbR$. Obviously,
$$
\pi_{x_0}(I-aR)=I-\pi_{x_0}(a)S=I-\pi_0(\hat{a})S=\pi_{0}(I-\hat{a}R).
$$
Consider the invertible isometry $J_{x_0}$ on $L_p(\bbR;E)$, defined
as  $(J_{x_0}f)(x)=f(x+x_0)$. Clearly,
$I-\hat{a}R=J_{x_0}(I-aR)J_{x_0}^{-1}$.
Hence, the operator $I-\hat{a}R$ is invertible if and only if
the operator $b=I-aR$ is invertible, and $\|(I-\hat{a}R)^{-1}\|=\|b^{-1}\|$.
Therefore, the estimate \eqref{binv} for $x=x_0$ follows
from the estimate \eqref{binv} for $x=0$.

Thus our purpose is to prove if $b=I-aR$ is invertible in $L_p(\bbR; E)$,
then the
operator $\pi_0 (b)=I-\pi_0(a)S$
is invertible in $l_p(\bbZ; E)$, and
\begin{equation}\label{est0}
\| [\pi_0 (b)]^{-1}\|_{B(l_p (\bbZ;E))} \leq \|b^{-1}\|_{B(L_p
(\bbR;E))}.
\end{equation}

We first show that
for any $(v_n) \in l_p (\bbZ;E)$ the
following estimate holds:
\begin{equation}
\| (v_n)\|_{l_p (\bbZ;E)} \leq \|b^{-1}\|_{B(L_p (\bbR;E))} \| \pi_0
(b) (v_n)\|_{l_p(\bbZ;E)}.
\lb{B}
\end{equation}
Let us fix a sequence
$(v_n)_{n\in\bbZ}\in l_p (\bbZ;E)$, a natural number $N>1$, and $\eps > 0$.

Recall that the function
$\bbR \ni x \mapsto a(x) v\in E$ is continuous for each $v\in
E$.  Choose $\del < 1$ such that
\begin{equation}
\| [a(x+n) -a(n) ] v_{n-1} \|_E < \eps,
\quad \forall x\in [0 , \del],\; n=-N, \ldots, N.
\lb{eps}
\end{equation}
Define $f\in L_p(\bbR;E)$ by $f(x) =v_n$
for $x\in [n, n+\del]$,
$|n| \leq N$, and $f(x) =0$ otherwise.  Since $b$ is an
invertible operator in $L_p
(\bbR;E)$, it follows that:
\begin{equation}
\|b^{-1} \|_{B(L_p)}^p \|bf\|_{L_p}^p \geq \|f\|_{L_p}^p =\sum_{n=-N}^N
\int_{n}^{n+\del} \|v_n\|^p \, dx = \del \sum_{n=-N}^N
\|v_n\|^p.
\lb{f}
\end{equation}
On the other hand, using \eqref{eps}, one has:
\begin{align*}
\|bf\|_{L_p}^p & =\int_{\bbR} \| f(x) -a(x) f(x-1)\|_E^p\, dx\\
&=\sum_{n=-N+1}^N \int_{n}^{n+\del} \|v_n -a(x) v_{n-1} \|_E^p \,
dx +\int_{-N}^{-N+\del} \| v_{-N}\|_E^p \, dx \\
&\quad + \int_{N+1}^{N+1+\del} \|a(x) v_N \|_E^p \, dx\\
& \leq \sum_{n=-N+1}^N \int_{0}^{\del}
\|v_n - a(n) v_{n-1} -[a(x+n)
-a(n)] v_{n-1} \|_E^p\, dx\\
&\quad +\del \|v_{-N}\|_E^p +\del \max_{x\in\bbR} \|a(x) \|^p \cdot
\|v_N\|_E^p\\
& \leq \del \sum_{n=-N}^N (\| v_n -a(n) v_{n-1} \|_E +\eps )^p +\del
\| v_{-N}\|_E^p +\del \max_{x\in\bbR} \|a(x) \|^p \|v_N\|_E^p.
\end{align*}
Combining this inequality with \eqref{f}, one has:
\begin{align*}
\sum_{n=-N}^N \|v_n\|^p & \leq \|b^{-1} \|_{B(L_p)}^p
\Big\{ \sum_{n=-N}^N
(\|v_n -a(n) v_{n-1}\|_E +\eps )^p \\
& \quad + \|v_{-N} \|_E^p +
\max_{x\in\bbR}\|a(x)\|^p \|v_N\|_E^p \Big\}.
\end{align*}
If $\eps \to 0$ and $N\to\infty$, then
$$
\|(v_n)\|_{l_p(\bbZ;E)} \leq \|b^{-1} \|_{B(L_p)} \Big(
\sum_{n=-\infty}^\infty \|v_n-a(n) v_{n-1} \|_E^p \Big)^{1/p}
$$
and \eqref{B} is proved.

Note, that \eqref{B} is sufficient to show \eqref{est0} provided
the operator $\pi_0(b)$ is invertible, and so it only remains
to show that $\pi_0 (b)$ is an operator onto $l_p(\bbZ;E)$.

Fix any $(v_n )\in l_p
(\bbZ;E)$ and let $f(x) =U(x,n-1)v_{n-1}$, $x\in [n-1/2, n+1/2)$,
$n\in\bbZ$.
Property~(iii) from the Definition~\ref{evfam} of
the evolutionary family
$\{U(x,s)\}$ implies that
$$\|U(x,n-1)\| \leq C e^{\beta
(x-n+1)} \leq Ce^{\frac{3}{2}\beta} \mbox{ for } x\in [n-1/2, n+1/2).$$
Hence
$f\in L_p(\bbR;E)$.
Since the operator
$b$, defined by $(bg)(x) =g(x) -U(x,x-1)g(x-1)$,
is invertible in $L_p (\bbR;E)$, it follows that there exists a unique
function $g\in L_p(\bbR;E)$ such that
\begin{equation}
g(x) -U(x,x-1) g(x-1)=f(x)
\lb{eq}
\end{equation}
for almost all $x\in\bbR$.  Since
$$
\|g\|_{L_p}^p =\sum_{n\in\bbZ} \int_{-n-1/2}^{-n+1/2} \|g(s)\|^p \, ds
=\int_{-1/2}^{1/2} \Big( \sum_{n\in\bbZ} \|g(s+n)\|^p\Big) \, ds <
\infty,
$$
the sequence $(g(s+n))_{n\in\bbZ}$ belongs to $l_p(\bbZ;E)$ for all
$s\in\Ome$ for some subset $\Ome \subset (-1/2, 1/2)$ of full measure.
For each $s\in\Ome$, let us define a function $h_s$ by the rule:
$$
h_s(x) =\begin{cases}
g(x), & \mbox{ if } n-\dfrac12 \leq x \leq n+s,\\
U(x,n+s) g(n+s), & \mbox{ if } n+s\leq x < n+\dfrac12,
\end{cases}
\quad n \in\bbZ.
$$
Clearly, $h_s \in L_p (\bbR;E)$ for each $s\in\Ome$ because
$(g(s+n))_{n\in\bbZ}\in l_p (\bbZ;E)$, and
$$\|U(x,n+s)\|
\leq Ce^{\beta (1/2 -s)}
\mbox{ for } x\in [n+s, n+1/2).$$

We note that $h_s$ is a solution of
equation~\eqref{eq}.  Indeed, for $x\in [n-1/2, n+s]$, equation \eqref{eq}
implies $h_s (x) -U(x,x-1)h_s(x-1) =f(x)$.  For $x\in [n+s, n+1/2)$, one has:
 \begin{align*}
& h_s(x) -U(x,x-1) h_s(x-1)  = \\
 & U(x,n+s) [g(n+s) -U(n+s, n-1+s)g(n-1+s)]\\
& = U(x,n+s) f(n+s) =U(x,n+s) U(n+s, n-1) v_{n-1}=f(x).
\end{align*}
But equation~\eqref{eq} has only one solution $g$ in $L_p (\bbR;E)$.
Hence $g=h_s$ for all $s\in\Ome$.

Fix $s<0$, $s\in\Ome$.  The function
$h_s(\cdot)$ is defined for $x=n$, $n\in\bbZ$.  Moreover, the sequence
$(h_s(n))_{n\in\bbZ}$ belongs to $l_p (\bbZ;E)$. Indeed, $h_s(n)
=U(n,n+s)g(n+s)$, the sequence $(g(n+s))_{n\in\bbZ} \in l_p (\bbZ;E)$,
and $\|U(n,n+s)\|\leq Ce^{-\beta s}$ by (iii) from Definition~\ref{evfam}.

Set $u_n =h_s (n) +v_n$.  Since $h_s(\cdot)$
satisfies the equation \eqref{eq} for $x=n$, $n\in \bbZ$, we have:
\begin{align*}
& \pi_0(b) (u_n)  = \\
& u_n -U(n,n-1) u_{n-1} =h_s (n) +v_n -U(n,n-1)
h_s(n-1)-U(n,n-1)v_{n-1} =\\
& v_n+ h_s(n) - U(n, n-1)h_s(n-1) - f(n)= (v_n).
\end{align*}
This identity proves that $\pi_0 (b)$ is an operator onto
$L_p(\bbZ;E)$.
\end{pf}

\begin{thm}\lb{cont}
Let $\{U(x,s)\}_{x\geq s}$ be an evolutionary family in a separable
Banach space $E$, and let
$\{e^{tG}\}_{t\geq 0}$ be the evolutionary semigroup acting on
$L_p(\bbR;E)$, $1\leq p < \infty$
by the rule $(e^{tG}f)(x) =U(x,x-t)f(x-t)$.  The
following conditions are equivalent:
\newcounter{ee}
\begin{list}{\arabic{ee})}{\usecounter{ee}}
\item $\{U(x,s)\}_{x\geq s}$ is (spectrally) hyperbolic in $E$;
\item $\sig (e^{tG}) \cap \bbT =\emp$, $t\neq 0$, in $L_p (\bbR;E)$;
\item $0\in\rho (G)$ in $L_p (\bbR;E)$.
\end{list}
Moreover, the Riesz projection $\calP$ that corresponds to the part $\sig
(e^G)\cap \bbD$ of the spectrum of the hyperbolic operator $e^G$ is related
to a strongly continuous, projection-valued function $P:\bbR \to B(E)$
that satisfies Definition~\ref{hyp} by the formula $(\calP f)(x) =P(x) f(x)$,
$x\in\bbR$, $f\in L_p(\bbR;E)$.
\end{thm}

\begin{pf}
2) $\Leftrightarrow$ 3) was proved in Theorem~\ref{t3.1}.

1) $\Rightarrow$ 2).
Without loss
of generality assume $t=1$.
>From the projection-valued function $P(\cdot)$ from
Definition~\ref{hyp}, let us define a
projection $\calP$ in $L_p(\bbR;E)$ by the rule $(\calP f)(x) =P(x) f(x)$.
Denote $\calQ =I-\calP$.  For $T=e^G$, condition a)
of Definition~\ref{hyp} implies $\calP T =T\calP$.
Set $T_P =\calP T \calP$ and $T_Q =\calQ T\calQ$. Then b) implies
$\sig (T_P)\subset
\bbD$ in $\Im \calP =\{f\in L_p (\bbR;E): f(x) \in \Im P(x)\}$.  Also b)
and c) imply that the operator $T_Q$,
which can be written as
$(T_Qf)(x) =Q(x) U(x,x-1) Q(x-1)f(x-1)$, is
an invertible operator, and
$\sig (T^{-1}_Q)\subset\bbD$ in $\Im \calQ =\Ker \calP$.
Hence, $\sig (T) \cap \bbT =\emp$.

2) $\Rightarrow$ 1).  Let $\frakB$ be
a Banach algebra with a norm $\|\cdot \|_1$
consisting of the operators $d$ on
$L_p(\bbR;E)$ of the form
$$
d=\sum_{k=-\infty}^\infty a_k R^k, \; a_k \in\frakA,\; \|d\|_1 :
=\sum_{k=-\infty}^\infty \| a_k\|_{B(L_p (\bbR;E))} < \infty.
$$

We first show that if $b=\lam -T$ is invertible in $L_p(\bbR;E)$ for all
$\lam \in \bbT$, then $(\lam -T)^{-1} \in \frakB$. This fact will be proved
in several steps.

First, without loss of generality let $\lam=1$.  Since $\sig (T) \cap
\bbT=\emp$, by Lemma~\ref{ling} the Riesz projection $\calP$ has a form
$(\calP f)(x)=P(x) f(x)$, where
the function $\bbR_0 \ni x \mapsto P(x) v\in E$ is a
bounded, measurable (in the strong sense in $E$) function for each $v\in E$,
where $\bbR_0$\ is a set of full measure in $\bbR$.
Recall also that $Q(x)=I-P(x)$.
Decompose $b=I-T=(\calP -T_P) \oplus (\calQ-T_Q)$.
Since $\sig (T_P) \subset
\bbD$ and $\sig (T^{-1}_Q )\subset \bbD$, one has that
$b^{-1}=(\calP -T_P)^{-1} \oplus (\calQ-T_Q)^{-1}$, where
\begin{equation}
(\calP -T_P )^{-1} =\sum_{k=0}^\infty T_P^k;
\quad (\calQ-T_Q)^{-1}=[-T_Q
(\calQ-T_Q^{-1})]^{-1} =-\sum_{k=-\infty}^{-1} T_Q^k.
\lb{neum}
\end{equation}
Notice
that $T_Q^{-1} =(\calQ aR\calQ)^{-1} =(QaQ(\cdot -1) R)^{-1} = R^{-1} (QaQ
(\cdot -1))^{-1} =[Q(\cdot +1) a(\cdot +1) Q(\cdot)]^{-1} R^{-1}$, and
that $T_P=aR\calP =aP(\cdot -1) R$.
Hence both operators $T_P^k$ and $T_Q^k$ can be written as $a_kR^k$
for some multiplication operators $a_k$.
The Neumann series in \eqref{neum} converge
absolutely.  Therefore,
\begin{equation}
b^{-1} =\sum_{k=-\infty}^\infty a_k R^k,\; \sum_{k=-\infty}^\infty
\|a_k\|_{B(L_p (\bbR;E))} < \infty,
\lb{invb}
\end{equation}
for each $k\in\bbZ$ the function $a_k: \bbR_0 \to B(E)$ is bounded, and
the function
$\bbR_0 \ni x\mapsto a_k(x) v\in E$ is measurable for each $v\in E$.

Our next aim is to show that the $a_k$ from \eqref{invb} belong to
$\frakA$, that is, the function
$x \mapsto a_k (x) v\in E$ extends to a continuous function from $\bbR$ for
each $v\in E$.  To this end let us define for $a_k$ from \eqref{invb} and
all $x\in\bbR_0$ the operator $\pi_x(a_k)$ in $l_p (\bbZ;E)$ as in
\eqref{pi}.  Denote:
\begin{equation}\label{pixb}
\pi_x(b^{-1})=\sum\pi_x(a_k) S^k,\quad
\pi_x(a_k) =\diag \{a_k (x+n)\}_{n\in\bbZ} \mbox{\quad for } x\in\bbR_0.
\end{equation}
Identities $bb^{-1}
=b^{-1} b=I$ in $L_p(\bbR;E)$ imply that $\pi_x (b)\cdot \pi_x
(b^{-1})=\pi_x(b^{-1}) \cdot \pi_x(b) =I$ in $l_p (\bbZ;E)$
for $x\in\bbR_0$.
Since the operator $b$ is invertible in $L_p(\bbR; E)$ by assumption,
 for each $x\in\bbR$ the operator $\pi_x(b)$ is
invertible in $l_p(\bbZ;E)$ by Lemma~\ref{inv}.  Hence
\begin{equation}
\pi_x(b^{-1}) =[\pi_x(b)]^{-1} \mbox{ for } x\in\bbR_0.
\lb{piinv}
\end{equation}

Recall that the function
$\bbR \ni x \mapsto a(x) v\in E$ is continuous for each $v\in
E$.  Also, the function $\bbR\ni x\mapsto \|a(x)\|\in \bbR_+$
is bounded. Hence for each
$(v_n) \in l_p (\bbZ;E)$ the function $\bbR \ni x \mapsto
\pi_x(b)(v_n) \in l_p (\bbZ;E)$ is continuous.  By Lemma~\ref{inv}
$\|[ \pi_x (b) ]^{-1}\|_{B(l_p)}$ are uniformly bounded for $x\in\bbR$.
This implies that the function
$\bbR \ni x \mapsto [\pi_x(b)]^{-1} (v_n) \in l_p (\bbZ;E)$
is continuous for each $(v_n ) \in l_p (\bbZ;E)$.
Indeed,
\begin{align*}
& \big|\big| \left([\pi_x(b)]^{-1}-[\pi_{x_0}(b)]^{-1}\right)
(v_n)\big|\big|_{l_p(\bbZ;E)} = \\
& \big|\big| [\pi_x(b)]^{-1}\cdot
\left[\pi_x(b)-\pi_{x_0}(b)\right]\cdot
[\pi_{x_0}(b)]^{-1}
(v_n)\big|\big|_{l_p(\bbZ;E)}
\end{align*}
for any $x,x_0 \in\bbR$, and the function
$\bbR \ni x \mapsto [\pi_x(b)]^{-1} (v_n) \in l_p (\bbZ;E)$
is continuous at $x=x_0$.

Fix $k_0 \in \bbZ$, $x_0 \in \bbR$, and $v\in E$.  Define $(\ti v_n) \in
l_p (\bbZ;E)$ as $\ti v_{-k_0} =v$ and $\ti v_n =0$ for $n\neq -k_0$.
Consider a sequence $x_m \to x_0$, $x_m \in\bbR_0$.  We will show that
$\{a_{k_0} (x_m) v\}_{m\in\bbN}$ is a Cauchy sequence in $E$ and will
define $a_{k_0} (x_0) v=\lim\limits_{m\to\infty} a_{k_0} (x_m) v$.  Then
the function
$\bbR \ni x \mapsto a_{k_0} (x) v\in E$ becomes a continuous function,
and $a_{k_0}\in\frakA$.

Note, that the $\pi_{x_m} (b^{-1})$ are defined
by the formula \eqref{pixb} since $x_m\in\bbR_0$.  For the
sequence $(\ti v_n)$ one has the following estimate:
\begin{align} \label{est5}
& \| [\pi_{x_{m'}} (b^{-1}) -\pi_{x_{m''}} (b^{-1})] (\ti v_n)
\|_{l_p}^p  \nonumber \\
&=\sum_{n\in\bbZ} \Big\|\sum_{k\in\bbZ} [a_k (x_{m'} +n) -a_k (x_{m''}+n)
]\ti v_{n-k} \Big\|_E^p \nonumber \\
& \geq \Big\| \sum_k [a_k (x_{m'}) -a(x_{m''})] \ti
v_{-k} \Big\|_E^p \nonumber \\
&= \| [a_{k_0} (x_{m'}) -a_{k_0} (x_{m''})] v\|_E^p.
\end{align}
Since $x_m \in\bbR_0$, formula~\eqref{piinv} is applicable.
Then the sequence
$$\pi_{x_m} (b^{-1}) (\ti v_n) =[\pi_{x_m} (b) ]^{-1} (\ti v_n)$$ is a
Cauchy sequence in $l_p (\bbZ;E)$ since the function
$\bbR\ni x \mapsto [\pi_x (b)
]^{-1}(\ti v_n) \in l_p (\bbZ;E)$ is continuous. In accordance with
\eqref{est5} the sequence  $\{a_{k_0} (x_m) v\}_{m\in\bbN}$ is a
Cauchy sequence in $E$, and $a_{k_0}\in\frakA$.

Since the $a_k$ from \eqref{invb} are continuous,
we have proved that $(\lam I -T)^{-1} \in
\frakB$ for all $\lam \in \bbT$.

The rest of the proof is standard
(cf.~\cite{2,18,27}).  Indeed, consider the absolutely convergent
Fourier series $f:\lam \mapsto \lam I -aR\lam ^0$ with the
coefficients from
$\frakB$.  For each $\lam \in\bbT$, the operator
$f(\lam) =b$ is invertible in
$\frakB$.  Hence the function $\bbT \ni \lam \mapsto [f(\lam)]^{-1}
\in\frakB$ is expandable (see, e.g., \cite{5}) into an absolutely
convergent Fourier series
$$
[f(\lam)]^{-1} =(\lam I -aR)^{-1} =\sum_{k=-\infty}^\infty d_k \lam^k,\;
\sum_k \|d_k\| < \infty, \; d_k \in\frakB.
$$
By the integral formula (see, e.g., \cite[p.~20]{9}) for
the Riesz projection $\calP$, we conclude that
$\calP =d_{-1}\in\frakB$.  Hence for some $a_k \in\frakA$ one has that
$$\calP=\sum_{k=-\infty}^{\infty} a_k R^k, \mbox { where }
\sum_{k=-\infty}^{\infty} \|a_k\|  <\infty.
$$

We will show that $a_k =0$ for $k\neq 0$, and $\calP
=a_0 \in \frakA$.  Indeed, by \eqref{com}, $\chi \calP =\calP \chi$ for
any bounded continuous scalar function $\chi$.  Then
$$
\chi\calP -\calP \chi =\sum_k a_k (\cdot ) [\chi (\cdot ) -\chi (\cdot
-k)] R^k =0.
$$
Then by picking $x_0$\ and $\chi$ such that $\chi (x_0) \neq \chi
(x_0 -k)$ for $k\ne0$, it follows that $a_k (x_0) =0$, $k\neq 0$.
\end{pf}

As we have mentioned above,
for the space $C_0(\bbR;E)$, part 1) $\Leftrightarrow$ 2) of
Theorem~\ref{cont} was proved in \cite{27}.

\end{document}